\documentclass[a4paper,11pt]{amsart}

\usepackage[utf8]{inputenc}
\usepackage{amsmath}
\usepackage{amssymb}
\usepackage{amsthm}
\usepackage{amsrefs}
\usepackage{mathrsfs}
\usepackage{xfrac}
\usepackage{pifont}

\usepackage{array,tabularx}
\usepackage{xcolor,tikz,tikz-cd}
\usetikzlibrary{arrows,arrows.meta}
\usetikzlibrary{calc}
\tikzcdset{arrow style=tikz, diagrams={>={Computer Modern Rightarrow[width=5.5pt,length=3pt]}}}
\usepackage{enumitem}
\usepackage{cleveref}
\usepackage[retainorgcmds]{IEEEtrantools}

\makeatletter
\def\thickhline{%
  \noalign{\ifnum0=`}\fi\hrule \@height \thickarrayrulewidth \futurelet
   \reserved@a\@xthickhline}
\def\@xthickhline{\ifx\reserved@a\thickhline
               \vskip\doublerulesep
               \vskip-\thickarrayrulewidth
             \fi
      \ifnum0=`{\fi}}
\makeatother

\theoremstyle{plain}
\newtheorem{theorem}{Theorem}[section]
\newtheorem*{theorem*}{Theorem}
\newtheorem{proposition}[theorem]{Proposition}
\newtheorem{lemma}[theorem]{Lemma}
\newtheorem{corollary}[theorem]{Corollary}

\theoremstyle{definition}
\newtheorem{definition}[theorem]{Definition}

\newtheorem{example}[theorem]{Example}
\newtheorem{remark}[theorem]{Remark}

  \def\ba{\begin{array}}
  \def\ea{\end{array}}
  \def\bc{\begin{corollary}}
  \def\ec{\end{corollary}}
  \def\bd{\begin{definition}}
  \def\ed{\end{definition}}
  \def\ben{\begin{enumerate}}
  \def\een{\end{enumerate}}
  \def\bse{\begin{equation*}}
  \def\ese{\end{equation*}}
  \def\be{\begin{example}}
  \def\ee{\end{example}}
  \def\bi{\begin{IEEEeqnarray*}}
  \def\ei{\end{IEEEeqnarray*}}
  \def\bit{\begin{itemize}}
  \def\eit{\end{itemize}}
  \def\bl{\begin{lemma}}
  \def\el{\end{lemma}}
  \def\bnn{\begin{notation}}
  \def\enn{\end{notation}}
  \def\bn{\begin{note}}
  \def\en{\end{note}}
  \def\bp{\begin{proposition}}
  \def\ep{\end{proposition}}
  \def\bq{\begin{proof}}
  \def\eq{\end{proof}}
  \def\br{\begin{remark}}
  \def\er{\end{remark}}
  \def\bs{\begin{solution}}
  \def\es{\end{solution}}
  \def\btab{\begin{table}}
  \def\etab{\end{table}}
  \def\btb{\begin{tabular}}
  \def\etb{\end{tabular}}
  \def\bt{\begin{theorem}}
  \def\et{\end{theorem}}

  \def\B{\mathrm{B}}

  \def\cG{\mathcal{G}}

  \def\cl{\colon}

  \def\id{\mathrm{id}}

  \def\Map{\mathrm{Map}}
  
  \def\N{\mathbb{N}}

  \def\PSp{\mathrm{PSp}}
  
  \def\PU{\mathrm{PU}}
  \def\p{\partial}
  \def\pr{\mathrm{pr}}

  \def\Sp{\mathrm{Sp}}
  \def\Sq{\mathrm{Sq}}
  \def\SU{\mathrm{SU}}

  \def\Spin{\mathrm{Spin}}
  
  \def\U{\mathrm{U}}

  \def\Z{\mathbb{Z}}

\DeclareMathOperator{\ev}{ev}

\DeclareMathOperator{\SO}{SO}

\DeclareFontFamily{U}{MnSymbolC}{}
\DeclareSymbolFont{MnSyC}{U}{MnSymbolC}{m}{n}
\DeclareMathSymbol{\diamondplus}{\mathbin}{MnSyC}{"7C}
\DeclareMathSymbol{\diamonddot}{\mathbin}{MnSyC}{"7E}
\DeclareFontShape{U}{MnSymbolC}{m}{n}{
    <-6>  MnSymbolC5
   <6-7>  MnSymbolC6
   <7-8>  MnSymbolC7
   <8-9>  MnSymbolC8
   <9-10> MnSymbolC9
  <10-12> MnSymbolC10
  <12->   MnSymbolC12}{}

\usepackage{xparse}
\makeatletter
\RenewDocumentCommand{\title}{om}{%
   \IfNoValueTF{#1}
     {\gdef\shorttitle{Homotopy types of $\Spin^c(n)$-gauge groups over $S^4$}}
     {\gdef\shorttitle{#1}}%
   \gdef\@title{#2}%
}
\makeatother

\title{Homotopy types of $\Spin^c(n)$-gauge groups over $S^4$}
\author{Simon Rea}
\date{\today}

\usepackage{libertine}
\usepackage[libertine]{newtxmath}
\renewcommand{\mathbb}{\varmathbb}

\begin{document}

\maketitle

\begin{abstract}
The gauge group of a principal $G$-bundle $P$ over a space $X$ is the group of $G$-equivariant homeomorphisms of $P$ that cover the identity on $X$. We consider the gauge groups of bundles over $S^4$ with $\Spin^c(n)$, the complex spin group, as structure group and show how the study of their homotopy types reduces to that of $\Spin(n)$-gauge groups over $S^4$. We then advance on what is known by providing a partial classification for $\Spin(7)$- and $\Spin(8)$-gauge groups over $S^4$.


\end{abstract}

\section{Introduction}

Let $G$ be a topological group and $X$ a space. The \emph{gauge group} $\cG(P)$ of a principal $G$-bundle $P$ over $X$ is defined as the group of $G$-equivariant bundle automorphisms of $P$ which cover the identity on $X$. A detailed introduction to the topology of gauge groups of bundles can be found in \cites{husemoller,piccinini1998conjugacy}. The
study of gauge groups is important for the classification of principal bundles, as well as understanding moduli spaces of connections on principal bundles \cites{cohen-milgram,theriault-moduli1,theriault-moduli2}.
% Gauge groups also play a key role in theoretical physics, where they are used to describe the parallel transport of point particles by means of connections on bundles. Famously, Donaldson \cite{donaldson} computed the rational cohomology of the classifying space of the gauge group of an $\SU(2)$-bundle over a simply-connected 4-manifold and used it to define a polynomial invariant of manifolds. This later allowed Donaldson to deduce strong results about the differential
% topology of 4-manifolds.
Gauge groups also play a key role in theoretical physics, where they are used to describe the parallel transport of point particles by means of connections on bundles. Famously, Donaldson \cite{donaldson} discovered a deep link between the gauge groups of certain $\SU(2)$-bundles and the differential topology of 4-manifolds.

Key properties of gauge groups are invariant under continuous deformation and so studying their homotopy theory is important. Having fixed a topological group $G$ and a space $X$, an interesting problem is that of classifying the possible homotopy types of the gauge groups $\cG(P)$ of principal $G$-bundles $P$ over $X$.

Crabb and Sutherland showed \cite[Theorem 1.1]{crabb-sutherland} that if $G$ is a compact, connected, Lie group and $X$ is a connected, finite CW-complex, then the number of distinct homotopy types of $\cG(P)$, as $P\to X$ ranges over all principal $G$-bundles over $X$, is finite. In fact, since isomorphic $G$-bundles give rise to homeomorphic gauge groups, it will suffice to the let $P\to X$ range over the set of isomorphism classes of principal $G$-bundles over $X$.

Explicit classification results have been obtained, especially for the case of gauge groups of bundles with low rank, compact, Lie groups as structure groups and $X=S^4$ as base space. 
In particular, the first such result was obtained by Kono \cite{kono91} in 1991. Using the fact that isomorphism classes of principal $\SU(2)$-bundles over $S^4$ are classified by $k\in\Z\cong\pi_3(\SU(2))$ and denoting by $\cG_k$ the gauge group of the principal $\SU(2)$-bundle $P_k\to S^4$ corresponding to the integer $k$, Kono showed that there is a homotopy equivalence $\cG_k\simeq\cG_l$ if, and only if, $(12,k)=(12,l)$, where $(m,n)$ denotes the greatest common divisor of $m$ and $n$. Since 12 has six divisors, it follows that there are precisely six homotopy types of $\SU(2)$-gauge groups over $S^4$. 

Results formally similar to that of Kono have been obtained for principal bundles over $S^4$ with different structure groups, among others, by: Hamanaka and Kono \cites{hamanaka-kono} for $\SU(3)$-gauge groups; Theriault \cite{theriault-su5,theriault-sun} for $\SU(n)$-gauge groups, as well as \cite{theriault-sp2} $\Sp(2)$-gauge groups;  Cutler \cite{cutler-sp3,cutler-unitary} for $\Sp(3)$-gauge groups and $\U(n)$-gauge groups; Kishimoto, Theriault and Tsutaya for $G_2$-gauge groups; Kamiyama, Kishimoto, Kono and Tsukuda for $\SO(3)$-gauge groups; Hasui, Kishimoto, Kono and Sato \cite{hasui16} for $\PU(3)$- and $\PSp(2)$-gauge groups; and Hasui, Kishimoto, So and Theriault \cite{hasui19} for bundles with exceptional Lie groups as structure groups. There are also several classification results for gauge groups of principal bundles with base spaces other than $S^4$ \cites{kono-tsukuda,claudio-spreafico,hamanaka-s6,kono-theriault,theriault-su3,membrillo-theriault,membrillo19,hamanaka-sp2,theriault-so,su(4)-s8,su(n)-s6,hasui16,rea21,mohammadi21,west17,huang-5d,huang-highd,so19-sun,so19}.

% Classification results for gauge groups of bundles with base space other than $S^4$ include: Kono and Tsukuda's \cite{kono-tsukuda} and Claudio and Spreafico's \cite{claudio-spreafico} extensions of \cite{kono91} to the cases of gauge groups of $\SU(2)$-bundles over simply connected 4-manifolds and over higher-dimensional spheres, respectively; Hamanaka and Kono \cites{hamanaka-s6} for $\SU(3)$-gauge groups over $S^6$; Kono and Theriault \cite{kono-theriault} for $\SU(3)$-gauge groups over $S^n$; Theriault \cite{theriault-su3} for $\SU(3)$-gauge groups over simply connected 4-manifolds; Hamanaka, Kaji and Kono \cite{hamanaka-sp2} for $\Sp(2)$-gauge groups over $S^8$; Theriault and So \cite{theriault-so} for $\Sp(2)$-gauge groups over closed, simply connected 4-manifolds; Mohammadi and Asadi-Golmankhaneh \cites{su(4)-s8,su(n)-s6} for $\SU(n)$-gauge groups over $S^6$ and $\SU(4)$-gauge groups over $S^8$. More recently, the results of \cite{hasui16} have been extended to bundles with higher, even dimensional spheres as base spaces by the author \cite{rea21} for $\PU(p)$-gauge groups and by Mohammadi \cite{mohammadi21} for $\PSp(2)$- and $\PSp(3)$-gauge groups.

The complex spin group $\Spin^c(n)$ was first introduced in 1964 in a paper of Atiyah, Bott and Shapiro \cite{atiyah-bott-shapiro}. There has been an increasing interest in the $\Spin^c(n)$ groups ever since the publication of the Seiberg-Witten equations for 4-manifolds \cite{witten94}, whose formulation requires the existence of $\Spin^c(n)$-structures, and more recently for the role they play in string theory \cites{bryant99,freed-witten,sati12}.

In this paper we examine $\Spin^c(n)$-gauge groups over $S^4$. We begin by recalling some basic properties of the complex spin group $\Spin^c(n)$ and showing that, provided $n\geq 3$, it can be expressed as a product of a circle and the real spin group $\Spin(n)$. For $n\geq 6$, we show that this decomposition is reflected in the corresponding gauge groups. 
\bt
\label{lem:ggdecomp}
For $n\geq 6$ and any $k\in\Z$, we have
\bse
\cG_k\bigl( \Spin^c(n)\bigr)\simeq S^1\times \cG_k\bigl(\Spin(n)\bigr).
\ese
\et
The homotopy theory of $\Spin^c(n)$-gauge groups over $S^4$ therefore reduced to that of the corresponding $\Spin(n)$-gauge groups. We advance on what is known on $\Spin(n)$-gauge groups by providing a partial classification for $\Spin(7)$- and $\Spin(8)$-gauge groups over $S^4$.

\bt
\label{thm:spin7}
(a) If $(168,k)=(168,l)$, there is a homotopy equivalence
\bse
\cG_k(\Spin(7))\simeq\cG_l(\Spin(7))
\ese 
after localising rationally or at any prime;

(b) If $\cG_k(\Spin(7))\simeq\cG_l(\Spin(7))$ then $(84,k)=(84,l)$.
\et

We note that the discrepancy by a factor of 2 between parts (a) and (b) is due to the same discrepancy for $G_2$-gauge groups. 

\bt
\label{thm:spin8}
(a) If $(168,k)=(168,l)$, there is a homotopy equivalence 
\bse
\cG_k(\Spin(8))\simeq\cG_l(\Spin(8))
\ese
after localising rationally or at any prime;

(b) If $\cG_k(\Spin(8))\simeq\cG_l(\Spin(8))$ then $(28,k)=(28,l)$. Furthermore, if $k$ and $l$ are multiples of 3, then $(3,k)=(3,l)$.
\et

For the $\Spin(8)$ case, in addition to the same 2-primary indeterminacy appearing in the $\Spin(7)$ case, there are also known \cites{kishimoto10,theriault-odd} difficulties at the prime 3 due to the non-vanishing of $\pi_{10}(\Spin(8))_{(3)}$.

\section{$\Spin^c(n)$ groups}

For $n\geq 1$, the complex spin group $\Spin^c(n)$ is defined as the quotient 
\bse
\frac{\Spin(n)\times \U(1)}{\Z/2\Z}
\ese
where $\Z/2\Z\cong\{(1,1),(-1,-1)\}\subseteq \Spin(n)\times \U(1)$ denotes the central subgroup of order 2. The group $\Spin^c(n)$ is special case of the more general notion of $\Spin^k(n)$ group introduced in \cite{albanese-spinh}.

The first low rank $\Spin^c(n)$ groups can be identified as follows:
\begin{itemize}
    \item $\Spin^c(1)\cong \U(1)\simeq S^1$;
    \item $\Spin^c(2)\cong \U(1)\times \U(1)\simeq S^1\times S^1$;
    \item $\Spin^c(3)\cong \U(2)\simeq S^1\times S^3$;
    \item $\Spin^c(4)\cong \{(A,B)\in \U(2)\times \U(2)\mid \det A= \det B\}$.
\end{itemize}
% \begin{table}[h!]
%     \centering
% \begin{tabular}{c|c}
%      $\ n\ $ & $\Spin^c(n)$ \\
%     \hline
%      1 & $\U(1)$ \\
%      2 & $\U(1)\times \U(1)$ \\
%      3 & $\U(2)$ \\
%      4 & $\{(A,B)\in \U(2)\times \U(2)\mid \det A= \det B\}$\\
% \end{tabular}
% \smallskip
% \caption{Low rank $\Spin^c(n)$ groups.}
% \end{table}
The group $\Spin^c(n)$ fits into a commutative diagram
\bse
\begin{tikzcd}
\{\pm 1\} \ar[d,hook]& \ar[l,"\pr_1"']\{(1,1),(-1,-1)\} \ar[r,"\pr_2"]\ar[d,hook] & \{\pm 1\}\ar[d,hook] \\
\Spin(n) \ar[d,"\lambda"]& \ar[l,"\pr_1"'] \Spin(n)\times S^1 \ar[r,"\pr_2"] \ar[d,"q"]& S^1 \ar[d,"2"] \\
\SO(n) & \ar[l]\Spin^c(n)\ar[r] & S^1,
\end{tikzcd}
\ese
where $q$ is the quotient map, $\lambda\cl \Spin(n)\to \SO(n)$ denotes the double covering map of the group $\SO(n)$ by $\Spin(n)$ and $2\cl S^1\to S^1$ denotes the degree 2 map. Furthermore, we observe that the map
\bse
\lambda\times 2\cl \Spin^c(n)\to \SO(n)\times S^1
\ese
is a double covering of $\SO(n)\times S^1$ by $\Spin^c(n)$.

\section{Method of classification}
\label{sec:method}

A principal bundle isomorphism determines a homeomorphism of gauge groups induced by conjugation \cite{piccinini1998conjugacy}. We therefore begin by considering isomorphism classes of principal $\Spin^c(n)$-bundles over $S^4$. These are classified by the free homotopy classes of maps $S^4\to \B\Spin^c(n)$. Since $\Spin^c(n)$ is connected, $\B\Spin^c(n)$ is simply-connected and hence there are isomorphisms 
\bse
[S^4,\B\Spin^c(n)]_{\text{free}} \cong \pi_3(\Spin^c(n))\cong\pi_3(\SO(n))\cong \begin{cases}0 & n=1,2\\ \Z^2 & n=4 \\ \Z & n=3, n\geq 5. \end{cases}
\ese

\br 
Note that for $n=3$ we have $\Spin^c(3)\cong \U(2)$, and the homotopy types of $\U(2)$-gauge groups over $S^4$ have been studied by Cutler in \cite{cutler-unitary}.
\er

For $n\geq 5$, let $\cG_k$ denote the gauge group of the  $\Spin^c(n)$-bundle $P_k\to S^4$ classified by $k\in\Z$.
By \cites{atiyah-bott,gottlieb}, there is a homotopy equivalence
\bse
\B \cG_k\simeq\Map_{k}(S^4,\B\Spin^c(n)),
\ese 
the latter space being the $k$-th component of $\Map(S^4,\B\Spin^c(n))$, meaning the connected component containing the map classifying $P_k\to S^4$.

There is an evaluation fibration
\bse
\Map^*_k(S^4,\B\Spin^c(n))\longrightarrow \Map_k(S^4,\B\Spin^c(n))\xrightarrow{\ \ev\ } \B\Spin^c(n),
\ese
where $\ev$ evaluates a map at the basepoint of $S^4$ and the fibre is the $k$-th component of the pointed mapping space $\Map^*(S^4,\B\Spin^c(n))$. This fibration extends to a homotopy fibration sequence
\bse
\cG_k\longrightarrow \Spin^c(n)\longrightarrow\Map^*_k(S^4,\B\Spin^c(n))\longrightarrow \B\cG_k\longrightarrow\B\Spin^c(n).
\ese
%where we used the equivalences $\B \cG_{i,k}\simeq\Map_{k}(S^{2i},\B\PU(n))$, $\Omega\B\PU(n)\simeq\PU(n)$, and $\Omega\B\cG_{i,k}\simeq\cG_{i,k}$.
Furthermore, by \cite{sutherland92} there is, for each $k\in\Z$, a homotopy equivalence 
\bse
\Map^*_k(S^4,\B\Spin^c(n))\simeq \Map^*_0(S^4,\B\Spin^c(n)).
\ese
The space on the right-hand side is homotopy equivalent to $\Map^*_0(S^3,\Spin^c(n))$ by the exponential law, and is more commonly denoted as $\Omega^3_0\Spin^c(n)$. We therefore have the following homotopy fibration sequence
\bse
\cG_k\longrightarrow \Spin^c(n)\xrightarrow{\ \p_k\ }\Omega_0^3\Spin^c(n)\longrightarrow \B\cG_k\longrightarrow\B\Spin^c(n),
\ese
which exhibits the gauge group $\cG_k$ as the homotopy fibre of the map $\p_k$. This is a key observation, as it implies that the homotopy theory of the gauge groups $\cG_k$ depends on the maps $\p_k$. 

\bl[Lang { \cite[Theorem 2.6]{lang-ehp}}]
\label{lem:lang}
The adjoint of $\p_k\cl \Spin^c(n)\to \Omega_0^3\Spin^c(n)$ is homotopic to the Samelson product $\langle k\epsilon,1\rangle\cl S^{3}\land\Spin^c(n)\to\Spin^c(n)$, where $\epsilon\in\pi_3(\Spin^c(n))$ is a generator and 1 denotes the identity map on $\Spin^c(n)$.\qed
\el

As the Samelson product is bilinear, we have $\langle k\epsilon,1\rangle \simeq  k\langle \epsilon,1\rangle$, and hence, taking adjoints once more, $\p_{k}\simeq k \p_{1}$.

\bl[Theriault {\cite[Lemma 3.1]{theriault-sp2}}]
\label{lem:theriault}
Let $X$ be a connected CW-complex and let $Y$ be an H-space with a homotopy inverse. Suppose that $f\in[X,Y]$ has finite order and let $m\in\N$ be such that $mf\simeq *$. Then, for any integers $k,l\in\Z$ such that $(m,k)=(m,l)$, the homotopy fibres of $k f$ and $lf$ are homotopy equivalent when localised rationally or at any prime. \qed
\el

\br
The lemma of Theriault is the local analogue of a lemma used by Hamanaka and Kono in their study \cite{hamanaka-kono} of $\SU(3)$-gauge groups over $S^4$.
\er

Part (a) of Theorems \ref{thm:spin7} and \ref{thm:spin8} will follow as applications of Lemma \ref{lem:theriault}. For parts (b) we will need to determine suitable homotopy invariants of the gauge groups.

\section{$\Spin^c(n)$-gauge groups}
We begin with a decomposition of $\Spin^c(n)$ as a product of spaces which will be reflected in an analogous decomposition of $\Spin^c(n)$-gauge groups.

\bl
\label{lem:s1fact}
For $n\geq 3$, we have $\Spin^c(n)\simeq S^1\times \widetilde \Spin^c(n)$, where $\widetilde \Spin^c(n)$ denotes the universal cover of $\Spin^c(n)$.
\el

\bq
We have $\pi_1(\Spin^c(n))\cong \Z$ for $n\geq 3$ (see, e.g.\ \cite{jost2008}).  By the Hurewicz and the universal coefficient theorems, we have isomorphisms
\bse
\Z\cong \pi_1(\Spin^c(n))\cong H_1(\Spin^c(n);\Z)\cong H^1(\Spin^c(n);\Z).
\ese
Therefore, we have maps $S^1\to \Spin^c(n)$ and $\Spin^c(n)\to K(\Z,1)\simeq S^1$ representing generators of $\pi_1(\Spin^c(n))$ and of $H^1(\Spin^c(n);\Z)$, respectively, such that the composite induces an isomorphism in $\pi_1$. Therefore, the homotopy fibration
\bse
\widetilde \Spin^c(n) \longrightarrow  \Spin^c(n) \longrightarrow K(\Z,1)\simeq S^1
\ese
defining the universal cover of $\Spin^c(n)$ admits a right homotopy splitting and hence, as $\Spin^c(n)$ is a group, we have
\bse\Spin^c(n)\simeq S^1\times \widetilde \Spin^c(n).\qedhere\ese
\eq

Note that as $\Spin^c(n)$ is a Lie group, we can equip $\widetilde\Spin^c(n)$ with a group structure for which the covering map
\bse
\varrho\cl \widetilde\Spin^c(n)\to \Spin^c(n)
\ese
is a group homomorphism.

\bl
\label{lem:univ}
For $n\geq 3$, we have $\widetilde \Spin^c(n)\simeq \Spin(n)$.
\el
\bq
Let $f\cl \Spin^c(n)\to K(\Z/2\Z,1)$ be the map in the fibration sequence
\begin{equation}
\label{eq:fibf}
\Z/2\Z\longrightarrow S^1\times \Spin(n)\xrightarrow{\ q \ } \Spin^c(n)\xrightarrow{\ f \ } K(\Z/2\Z,1) %\tag{$\star$}
\end{equation}
arising from the double covering of $\Spin^c(n)$ and let $g\cl \Spin^c(n)\to K(\Z,1)$ be the generator of $H^1(\Spin^c(n);\Z)$ realising the splitting of $\Spin^c(n)$ in Lemma \ref{lem:s1fact}.
We claim that there is a homotopy commutative diagram
\bse
\begin{tikzcd}
\Spin^c(n)\ar[d,"g"] \ar[r,"f"]& K(\Z/2\Z,1)\ar[d,equal]\\
K(\Z,1) \ar[r,"\rho"]& K(\Z/2\Z,1)
\end{tikzcd}
\ese
where $\rho$ is the map induced by the mod 2 reduction $\Z\to\Z/2\Z$. Since the target space is an Eilenberg-MacLane space, it will be sufficient to check cohomology.

Indeed, on the one hand, since $g$ represents the generator of $H^1(\Spin^c(n); \Z)$, the composite $\rho\circ g$ represents the unique class in $H^1(\Spin^c(n);\Z/2\Z)\cong \Z/2\Z$ (cf.\ Harada and Kono \cite{konoBSpinc} for the mod 2 cohomology of $\Spin^c(n)$).

On the other hand, applying $\pi_1$ to the fibration sequence \eqref{eq:fibf} we obtain an exact sequence
\bse
\pi_1(\Z/2\Z)\longrightarrow \pi_1( S^1\times \Spin(n)) \xrightarrow{\ q_* \ } \pi_1( \Spin^c(n))\xrightarrow{\ f_* \ } \pi_1(K(\Z/2\Z,1)).%\longrightarrow\pi_1(\B S^1\times \B \Spin(n)) 
\ese
Recalling from Section 2 that $q_*$ induces multiplication by 2 on the fundamental groups, we therefore have
\bse
0 \longrightarrow \Z\xrightarrow{\  q_*\,=\,2 \ } \Z \xrightarrow{\ f_* \ } \Z/2\Z\longrightarrow 0 
\ese
since $\pi_1(\B S^1\times \B \Spin(n))\cong 0$. Hence $f_*$ is reduction mod 2 on $\pi_1$.
Applying the Hurewicz theorem and changing coefficients to $\Z/2\Z$ then gives a commutative diagram
\bse
\begin{tikzcd}
\pi_1( \Spin^c(n))\ar[r,"\pi_1(f)"] \ar[d,"h_1","\cong"']&\pi_1(K(\Z/2\Z,1))\ar[d,"h_1","\cong"']\\
H_1(\Spin^c(n);\Z)\ar[r,"H_1(f)"] \ar[d] & H_1(K(\Z/2\Z,1);\Z)\ar[d,"\cong"']\\
H_1(\Spin^c(n);\Z/2\Z)\ar[r,"f_*"] & H_1(K(\Z/2\Z,1);\Z/2\Z).
\end{tikzcd}
\ese
Since $\pi_1(f)$ and the composite in the left column are both reduction mod 2, the diagram implies that $H_1(f)$ is also reduction mod 2. Hence $f_*$ is an isomorphism in mod 2 homology. Finally, by the universal coefficient theorem for cohomology with field coefficients, we see that 
\bse
f^*\cl H^1(K(\Z/2\Z,1);\Z/2\Z)\longrightarrow H^1(\Spin^c(n);\Z/2\Z)
\ese
is an isomorphism. Therefore  $f\cl \Spin^c(n)\to K(\Z/2\Z,1)$ also represents the unique class in $H^1(\Spin^c(n);\Z/2\Z)\cong \Z/2\Z$ and hence we have $f\simeq \rho\circ g$. Taking fibres, we obtain a diagram of homotopy fibrations
\bse
\begin{tikzcd}
\widetilde  \Spin^c(n) \ar[r,equal]\ar[d,"\psi"]&\widetilde  \Spin^c(n) \ar[d,"\varrho"]\ar[r]&*\ar[d]\\
S^1 \times \Spin(n) \ar[r,"q"] \ar[d,]& \Spin^c(n)\ar[d,"g"] \ar[r,"f"]& K(\Z/2\Z,1)\ar[d,equal]\\
S^1 \ar[r,"2"]& S^1\simeq K(\Z,1) \ar[r,"\rho"]& K(\Z/2\Z,1)
\end{tikzcd}
\ese
which defines a map $\psi\cl\widetilde  \Spin^c(n)\to S^1 \times \Spin(n)$. In particular, the fibration in the leftmost column induces an exact sequence
\bse
\pi_m(\Omega S^1)\longrightarrow\pi_m(\widetilde  \Spin^c(n)) \xrightarrow{\ \psi_*\ } \pi_m(S^1 \times \Spin(n)) \longrightarrow \pi_m(S^1)
\ese
for each $m>1$. Given that the projection $\pr_2\cl S^1\times \Spin(n)\to \Spin(n)$ induces isomorphisms $\pi_m(S^1 \times \Spin(n))\cong\pi_m(\Spin(n))$ for $m>1$ and that the groups $\widetilde  \Spin^c(n)$ and $\Spin(n)$ are both simply-connected, the composite $\pr_2\circ\psi$ induces isomorphisms on all homotopy groups and is therefore a homotopy equivalence by Whitehead's theorem. Hence $\widetilde \Spin^c(n)\simeq \Spin(n)$.
\eq

We are now ready to show that the decomposition 
\bse
\cG_k(\Spin^c(n))\simeq S^1\times\cG_k(\Spin(n)) \ese
for $n\geq 6$ holds as stated in Theorem \ref{lem:ggdecomp}.

\bq[Proof of Theorem \ref{lem:ggdecomp}] 
Identifying the universal cover of $\Spin^c(n)$ as $\Spin(n)$ as in Lemma \ref{lem:univ}, there is a covering fibration 
\bse 
\Spin(n) \xrightarrow{\ \varrho \ }\Spin^c(n) \xrightarrow{\ g \ } S^1
\ese
where $\varrho$ is a group homomorphism. Let $s\cl S^1 \to \Spin^c(n)$ be a right homotopy inverse of $g$, which exists by Lemma \ref{lem:s1fact}.

As $\pi_4(\Spin^c(n))\cong 0$ for $n\geq 6$, there is a lift in the diagram
\bse
\begin{tikzcd}
 & \ar[dl,dashed,"a"'] S^1 \ar[d,"s"]\ar[dr,"*"]&\\
 \cG_k(\Spin^c(n)) \ar[r]& \Spin^c(n) \ar[r,"\p_k"] & \Omega^3_0\Spin^c(n) .
\end{tikzcd}
\ese
% \bse
% \begin{tikzcd}
% &\cG_k(\Spin^c(n)) \ar[d] \\
% \ar[ur,dashed,"a"'] S^1 \ar[r,hook]\ar[dr,"*"]&  \Spin^c(n) \ar[d,"\p_k"] \\
% & \Omega^3_0\Spin^c(n) .
% \end{tikzcd}
% \ese

% --------- 12 july 2021

% Therefore, by \Cref{lem:s1fact,lem:univ}, letting $b$ be the map in the diagram
% \bse
% \begin{tikzcd}[row sep=large]
%  & \ar[dl,"a"'] S^1 \ar[dd,bend left=80,"{\id}"]\ar[d,"s"]\\
%  \cG_k(\Spin^c(n)) \ar[dr,"b"]\ar[r]&\Spin^c(n)\ar[d,"g"]\\
%  & S^1,
% \end{tikzcd}
% \ese

% we have $\cG_k(\Spin^c(n))\simeq S^1\times F_b$, where $F_b$ denotes the homotopy fibre of $b$. 
%--------

Define the map $b$ to be the composite
\bse
\cG_k(\Spin^c(n)) \longrightarrow \Spin^c(n) \xrightarrow{\ g \ } S^1.
\ese
Since $s$ is a right homotopy inverse for $g$, the map $a$ is a right homotopy inverse for $b$. Therefore we have $\cG_k(\Spin^c(n))\simeq S^1\times F_b$, where $F_b$ denotes the homotopy fibre of $b$.
% \bse
% \begin{tikzcd}
% F \ar[d]\ar[r]& \ar[d] \Spin(n) \ar[r,"\tilde \p_k"] & \Omega^3_0\widetilde\Spin^c(n)\ar[d,"\simeq"]\\
%  \cG_k(\Spin^c(n)) \ar[r]& \Spin^c(n)\ar[r,"\p_k"] &\Omega^3_0\Spin^c(n)
% \end{tikzcd}
% \ese

As the covering map $\varrho\cl \Spin(n)\to\Spin^c(n)$ is a group homomorphism, it classifies to a map
\bse
\B \varrho\cl \B \Spin(n) \to \B\Spin^c(n).
\ese
Since $\varrho$ induces an isomorphism in $\pi_3$, it respects path-components in $\Map_k(S^4,-)$ and $\Map^*_k(S^4,-)$ for any $k\in\Z$. We therefore have a diagram of fibration sequences
\begin{equation}
\label{eq:fib-diag}
\begin{tikzcd}[column sep=14pt]
\cdots\ar[r]&\Map^*_k(S^4,\B\Spin(n))\ar[r]\ar[d,"(\B \varrho)_*"] & \Map_k(S^4,\B\Spin(n))\ar[d,"(\B \varrho)_*"]\ar[r]& \B\Spin(n)\ar[d,"\B \varrho"]\\
\cdots\ar[r]&\Map^*_k(S^4,\B\Spin^c(n))\ar[r] & \Map_k(S^4,\B\Spin^c(n))\ar[r]&\B \Spin^c(n).
\end{tikzcd}
\end{equation}
Furthermore, observe that for all $k\in\Z$ we have
\bse
\pi_m(\Map^*_k(S^4,\B\Spin(n)))\cong \pi_m(\Omega_0^3\Spin(n))\cong\pi_{m+3}(\Spin(n))
\ese
and, similarly, $\pi_m(\Map^*_k(S^4,\B\Spin^c(n)))\cong\pi_{m+3}(\Spin^c(n))$. Since $\varrho$ induces isomorphisms on $\pi_m$ for $m\geq 2$, it follows that $(\B\varrho)_*$ induces isomorphisms
\bse
\pi_m((\B\varrho)_*)\cl \pi_m(\Map^*_k(S^4,\B\Spin(n))) \xrightarrow{ \ \cong \ }\pi_m(\Map^*_k(S^4,\B\Spin^c(n)))
\ese
for all $m$ and is therefore a homotopy equivalence by Whitehead's theorem.

We can extend the fibration diagram \eqref{eq:fib-diag} to the left as
\bse
\begin{tikzcd}
\cG_k(\Spin(n))\ar[d,"\cG_k(\varrho)"]\ar[r] & \Spin(n)\ar[d,"\varrho"]\ar[r,"\p_k'"] & \Map^*_k(S^4,\B\Spin(n))\ar[d,"(\B \varrho)_*","\simeq"']\ar[r] & \cdots\\
\cG_k(\Spin^c(n))\ar[r] & \Spin^c(n)\ar[r,"\p_k"] & \Map^*_k(S^4,\B\Spin^c(n))\ar[r] & \cdots
\end{tikzcd}
 \ese
% \bse
% \begin{tikzcd}%[column sep=18pt]
% \cG_k(\Spin(n))\ar[d,"\cG_k(\varrho)"]\ar[r]&\Spin(n)\ar[d,"\varrho"]\ar[r,"\p_k'"]&\Map^*_k(S^4,\B\Spin(n))\ar[d,"(\B \varrho)_*","\simeq"']\ar[r]&\cdots\\
% \cG_k(\Spin^c(n))\ar[r]\ar[d,"b"]&\Spin^c(n)\ar[d,"g"]\ar[r,"\p_k"]&\Map^*_k(S^4,\B\Spin^c(n))\ar[r]\ar[d]&\cdots\\
% S^1 \ar[r,equal]& S^1\ar[r] &*, &
% \end{tikzcd}
% \ese
where $\p_k'$ denotes the boundary map associated to $\Spin(n)$-gauge groups over $S^4$.

Since $(\B \varrho)_*$ is a homotopy equivalence, the leftmost square is a homotopy pull-back. Since we know that there is a fibration
\bse 
\Spin(n) \xrightarrow{\ \varrho \ }\Spin^c(n) \xrightarrow{\ g \ } S^1,
\ese
it follows that we also have a fibration
\bse 
\cG_k(\Spin(n)) \xrightarrow{\ \cG_k(\varrho) \ }\cG_k(\Spin^c(n)) \xrightarrow{\ b \ } S^1.
\ese
In particular, the space $\cG_k(\Spin(n))$ is seen to be the homotopy fibre $F_b$ of the map $b\cl \cG_k(\Spin^c(n))\to S^1$ and hence we have
\bse
\cG_k(\Spin^c(n))\simeq S^1\times\cG_k(\Spin(n)). \qedhere
\ese
\eq
In light of Theorem \ref{lem:ggdecomp}, the homotopy theory of $\Spin^c(n)$-gauge groups over $S^4$ for $n\geq 6$ is completely determined by that of $\Spin(n)$-gauge groups over $S^4$.

\br
By a result of Cutler \cite{cutler-unitary}, there is a decomposition 
\bse
\cG_k(\U(2))\simeq S^1\times \cG_k(\SU(2))
\ese
of $\U(2)$-gauge groups over $S^4$ whenever $k$ is even. Given that $\Spin^c(3)\cong \U(2)$ and $\Spin(3)\cong \SU(2)$, the statement of Theorem \ref{lem:ggdecomp} still holds true when $n=2$ provided that $k$ is even. Cutler also shows that $\cG_k(\U(2))\simeq S^1\times \cG_k(\PU(2))$ for odd $k$, so Theorem \ref{lem:ggdecomp} does not hold for $n=2$.
\er

\section{$\Spin(n)$-gauge groups}

We now shift our focus to principal $\Spin(n)$-bundles over $S^4$ and the classification of their gauge groups. In the interest of completeness, we recall that, for $n\leq 6$, the following exceptional isomorphisms hold.
% \begin{table}[h!]
%     \centering
% \begin{tabular}{c|cccccc}
%     $n$ &   1 & 2 & 3 & 4 & 5 & 6\\
%     \hline
%     $\Spin(n)$ & $\mathrm{O}(1)$ & $\U(1)$  & $\SU(2)$ & $\SU(2)\times \SU(2)$ & $\Sp(2)$ & $\SU(4)$
% \end{tabular}
% \caption{The exceptional isomorphisms.}
% \end{table}
\begin{table}[h!]
    \centering
\begin{tabular}{c|c}
     $\ n\ $ & $\Spin(n)$ \\
    \hline
     1 & $\mathrm{O}(1)$ \\
     2 & $\U(1)$ \\
     3 & $\SU(2)$ \\
     4 & $\SU(2)\times \SU(2)$\\
     5 & $\Sp(2)$\\
     6 & $\SU(4)$
\end{tabular}
\smallskip
\caption{The exceptional isomorphisms.}
\end{table}

The cases $n=1,2$ are trivial. Indeed, as $\pi_3(\mathrm{O}(1))\cong\pi_3(\U(1))\cong 0$, there is only one isomorphism class of $\mathrm{O}(1)$- and $\U(1)$-bundles over $S^4$ (namely, that of the trivial bundle), and hence there is only one possible homotopy type for the corresponding gauge groups. The case $n=3$ was studied by Kono in \cite{kono91}. The case $n=4$ can be reduced to the $n=3$ case by \cite[Theorem 5]{booth-remarks}. The case $n=5$ was studied by Theriault in \cite{theriault-sp2}. Finally, the case $n=6$ was studied by Cutler and Theriault in \cite{cutler-theriault}.

We shall now explore the $n=7$ case. Recall that we have a fibration sequence
\bse
\cG_k(\Spin(7))\longrightarrow \Spin(7)\xrightarrow{\ k\p_1\ } \Omega^3_0\Spin(7).
\ese

\bl
\label{lem:spin7pnot2}
Localised away from the prime 2, the boundary map 
\bse
\Spin(7)\xrightarrow{\ \p_1\ } \Omega^3_0\Spin(7)
\ese
has order $21$.
\el

\bq
Harris \cite{harris61} showed that $\Spin(2m+1)\simeq_{(p)}\Sp(m)$ for odd primes $p$. This result was later improved by Friedlander \cite{friedlander} to a $p$-local homotopy equivalence of the corresponding classifying spaces. Then, in particular, localising at an odd prime $p$, we have a commutative diagram
\bse
\begin{tikzcd}
\Spin(7)\ar[r,"\p_1"]\ar[d,"\simeq"] & \Omega_0^3\Spin(7)\ar[d,"\simeq"]\ar[r] & \Map_1(S^4,\B\Spin(7))\ar[d,"\simeq"] \ar[r]& \B \Spin(7)\ar[d,"\simeq"]\\
\Sp(3)\ar[r,"\p_1'"] & \Omega_0^3\Sp(3) \ar[r]&\Map_1(S^4,\B\Sp(3)) \ar[r]& \B \Sp(3)
\end{tikzcd}
\ese 
where $\p_1'\cl\Sp(3)\to\Omega_0^3\Sp(3)$ denotes the boundary map associated to $\Sp(3)$-gauge groups over $S^4$ studied in \cite{cutler-sp3}. Hence the result follows from the calculation in \cite[Theorem 1.2]{cutler-sp3} where it is shown that $\p_1'$ has order 21 after localising away from the prime 2.
\eq

\bl
\label{lem:facto}
Let $F\to X\to Y$ be a homotopy fibration, where $F$ is an H-space, and let $\p\cl \Omega Y \to F$ be the homotopy fibration connecting map. Let $\alpha\cl A\to \Omega Y$ and $\beta\cl B\to \Omega Y$ be maps such that
\ben
\item $\mu \circ (\alpha\times \beta)\cl A\times B\to \Omega Y$ is a homotopy equivalence, where $\mu$ is the loop multiplication on $\Omega Y$;
\item $\p\circ \beta\cl B\to \Omega Y$ is nullhomotopic.
\een
Then the orders of $\p$ and $\p\circ \alpha$ coincide.
\el

\bq
Let $\theta\cl \Omega Y\times F\to F$ denote the canonical homotopy action of the loopspace $\Omega Y$ onto the homotopy fibre $F$, and let $e=\mu\circ (\alpha\times \beta)$. Consider the diagram
\bse
\begin{tikzcd}[row sep=large]
& A\times B \ar[dl,"\pr_1"']\ar[d,"\alpha\times\beta"] \ar[dr,"e"]&\\
A \ar[d,"\alpha"]& \Omega Y\times \Omega Y\ar[d,"{\id}\times\p"]\ar[r,"\mu"] & \Omega Y\ar[d,"\p"]\\
\Omega Y \ar[r,hook] \ar[rr,bend right=35,"\p"]& \Omega Y \times F \ar[r,"\theta"] & F.
\end{tikzcd}
\ese
The left portion of the diagram commutes by the assumption that $\p\circ\beta\simeq *$, while the right and bottom portions commute by properties of the canonical action $\theta$. Therefore
\bse
\p\simeq\p\circ\alpha\circ\pr_1\circ e^{-1},
\ese
and hence the orders of $\p$ and $\p\circ \alpha$ coincide.
\eq

\bl
\label{lem:spin7p=2}
Localised at the prime 2, the order of the boundary map 
\bse
\Spin(7)\xrightarrow{\ \p_1\ } \Omega^3_0\Spin(7)
\ese
is at most $8$.
\el

\bq
The strategy here will be to show that $\p_8$ is nullhomotopic. This will suffice as we have $\p_8\simeq8\p_1$ by Lemma \ref{lem:lang}.

By a result of Mimura \cite{mimura67}*{Proposition 9.1}, the fibration
\bse
G_2\xrightarrow{\ \alpha \ } \Spin(7) \longrightarrow S^7
\ese
splits at the prime 2. Let $\beta\cl S^7\to \Spin(7)$ denote a right homotopy inverse for $\Spin(7)\to S^7$. Then the composite
\bse
G_2\times S^7 \xrightarrow{\ \alpha\, \times \, \beta \ }\Spin(7)\times\Spin(7)\xrightarrow{\ \mu \ }\Spin(7)
\ese
is a 2-local homotopy equivalence.

Observe that we have $\p_8\circ\beta\simeq*$ since $\pi_{10}(\Spin(7))\cong\Z/8\Z$ and $\p_8\circ\beta\simeq8\p_1\circ\beta$. Therefore, by \Cref{lem:facto}, the order of $\p_8$ equals the order of $\p_8 \circ \alpha$.
% \bse
% \begin{tikzcd}[column sep=large]
% G_2 \times S^7 \ar[d,"\pr_1'"]\ar[r,"{\mu\circ(\alpha\times\beta)}"]& \Spin(7)\ar[r,"\p_8"] & \Omega^3_0\Spin(7)\ar[d,equal]\\
% G_2\ar[rr,"\p_8\circ\alpha"]&&  \Omega^3_0\Spin(7).
% \end{tikzcd}
% \ese
As $\alpha$ is a group homomorphism, there is a diagram of evaluation fibrations
\bse
\begin{tikzcd}
G_2 \ar[r,"\p_8'"]\ar[d,"\alpha"]&\Omega^3_0G_2 \ar[d,"\Omega^3\alpha"]\\
\Spin(7)\ar[r,"\p_8"] & \Omega^3_0\Spin(7).
\end{tikzcd}
\ese
Since $\p_8'\simeq 8\p'_1\simeq *$ by \cite[Theorem 1.1]{g_2}, we must have $\p_8\simeq *$.
\eq

\bq[Proof of Theorem \ref{thm:spin7} (a)]
Lemmas \ref{lem:spin7pnot2} and \ref{lem:spin7p=2} imply that $168\p_1\simeq *$, so the result follows from Lemma \ref{lem:theriault}.
\eq

We now move on to consider $\Spin(8)$-gauge groups. 
\bl
\label{lem:spin8p2}
Localised at the prime 2 (resp. 3), the order of the boundary map 
\bse
\Spin(8)\xrightarrow{\ \p_1\ } \Omega^3_0\Spin(8)
\ese
is at most $8$ (resp. 3).
\el

\bq
There is a fibration
\bse
\Spin(7)\longrightarrow \Spin(8)\longrightarrow S^7
\ese
which splits after localisation at any prime. Therefore, we have a homotopy equivalence $\Spin(8)\simeq \Spin(7)\times S^7$ realised by maps $\alpha\cl\Spin(7)\to\Spin(8)$ and $\beta\cl S^7\to\Spin(8)$, where $\alpha$ is a group homomorphism. Integrally, we have $\pi_{10}(\Spin(8))\cong\Z/24\Z\oplus\Z/8\Z$ (see, e.g.\ the table in \cite{mimura95}). Hence the same argument presented in the proof of Lemma \ref{lem:spin7p=2} shows that $8\p_1\simeq *$ and $3\p_1\simeq *$ after localising at $p=2$ and $p=3$, respectively.
\eq

\bl
\label{lem:spin8pnot2}
Let $p\neq 3$ be an odd prime. Then the $p$-primary orders of the maps $\p_1\cl \Spin(7)\to\Omega^3_0\Spin(7)$ and $\p_1\cl \Spin(8)\to\Omega^3_0\Spin(8)$ coincide.
\el

\bq
As $\pi_{10}(\Spin(8))\cong\Z/24\Z\oplus\Z/8\Z$, any map $S^7 \to \Spin(8)$ is nullhomotopic after localisation at an odd prime $p$ different from 3. Thus, decomposing $\Spin(8)$ as $\Spin(7)\times S^7$ and arguing as in the proof of Lemma \ref{lem:spin7p=2} yields the statement.
\eq

\bq[Proof of Theorem \ref{thm:spin8} (a)]
Lemmas \ref{lem:spin8p2} and \ref{lem:spin8pnot2} imply that $168\p_1\simeq *$, so the result follows from Lemma \ref{lem:theriault}.
\eq

\section{Homotopy invariants of $\Spin(n)$-gauge groups}

\bl
\label{lem:spin7b1}
If $\cG_k(\Spin(7))\simeq \cG_l(\Spin(7))$, then $(21,k)=(21,l)$.
\el

\bq
As in the proof of Lemma \ref{lem:spin7pnot2}, localising at an odd prime, we have an equivalence $\B \Spin(7)\simeq_{(p)}\B \Sp(3)$. We therefore have a diagram of homotopy fibrations
\bse
\begin{tikzcd}
\Spin(7)\ar[r,"\p_k"]\ar[d,"\simeq"] & \Omega_0^3\Spin(7)\ar[d,"\simeq"] \ar[r]& \B\cG_k(\Spin(7))\ar[d]\ar[r] &\B\Spin(7)\ar[d,"\simeq"]\\
\Sp(3)\ar[r,"\p_k'"] & \Omega_0^3\Sp(3)\ar[r] & \B\cG_k(\Sp(3)) \ar[r]&\B\Sp(3)
\end{tikzcd}
\ese 
where $\p_k'\cl\Sp(3)\to\Omega_0^3\Sp(3)$ denotes the boundary map studied in \cite{cutler-sp3}. Thus, by the five lemma, we have
\bse
\pi_{11}(\B\cG_k(\Spin(7)))\cong \pi_{11}(\B\cG_k(\Sp(3))).
\ese
Hence the result now follows from the calculations in \cite[Theorem 1.1]{cutler-sp3} 
%\cite{cutler-sp3,bott56
where it is shown that, integrally,
\bse
\pi_{11}(\B\cG_k(\Sp(3)))\cong \Z/(84,k)\Z. \qedhere
\ese
\eq

In their study of the homotopy types of $G_2$-gauge groups over $S^4$ in \cite{g_2}, Kishimoto, Theriault and Tsutaya constructed a space $C_k$ for which 
\bse
H^*(C_k)\cong H^*(\cG_k(G_2))
\ese
in mod 2 cohomology in dimensions 1 through 6. The cohomology of $C_k$ is then shown to be as follows.

\bl[\cite{g_2}*{Lemma 8.3}]
\label{lem:c_k}
We have
\begin{itemize}
    \item if $(4,k)=1$ then $C_k\simeq S^3$, so $H^*(C_k)\cong H^*(S^3)$;
    \item if $(4,k)=2$ or $(4,k)=4$ then $    H^*(C_k)\cong H^*(S^3)\oplus H^*(P^5(2))\oplus H^*(P^6(2))$, where $P^n(p)$ denotes the $n$th dimensional mod $p$ Moore space;
\item if $(4,k)=2$ then $\Sq^2$ is non-trivial on the degree 4 generator in $H^*(C_k)$;
\item if $(4,k)=4$ then $\Sq^2$ is trivial on the degree 4 generator in $H^*(C_k)$. \hfill$\square$
\end{itemize}
\el
We make use of the same spaces $C_k$ as follows.
\bl
\label{lem:spin7b2}
If $\cG_k(\Spin(7))\simeq\cG_l(\Spin(7))$, then we have $(4,k)=(4,l)$.
\el
\bq
As in the proof of Lemma \ref{lem:spin7p=2}, recall that we have a 2-local homotopy equivalence
\bse
G_2\times S^7 \xrightarrow{\ \alpha\, \times \, \beta \ }\Spin(7)\times\Spin(7)\xrightarrow{\ \mu \ }\Spin(7).
\ese
Since the map $\alpha\cl G_2\to\Spin(7)$ is a homomorphism, we have a commutative diagram
\bse
\begin{tikzcd}[column sep=large]
G_2 \ar[d,"\alpha"]\ar[r,"\p_1'"]& \Omega^3_0G_2\ar[d,"\Omega^3\alpha"]\\
\Spin(7)\ar[r,"\p_1"] & \Omega^3_0\Spin(7).
\end{tikzcd}
\ese
Furthermore, as $\pi_7(\Omega_0^3G_2)\cong\pi_{10}(G_2)\cong0$, we have
\bse
\pi_7(\Omega^3_0\Spin(7))\cong \pi_7(\Omega_0^3G_2)\oplus \pi_7(\Omega^3S^7)\cong \pi_7(\Omega^3S^7),
\ese
and thus there is a commutative diagram
\bse
\begin{tikzcd}[column sep=large]
S^7 \ar[d,"\beta"]\ar[r,"\gamma"]& \Omega^3S^7\ar[d,"\Omega^3\beta"]\\
\Spin(7)\ar[r,"\p_1"] & \Omega^3_0\Spin(7)
\end{tikzcd}
\ese
for some $\gamma$ representing a class in $\pi_7(\Omega^3S^7)\cong\pi_{10}(S^7)\cong\Z/8\Z$.

% Since the lowest dimensional cell in $G_2\times S^7/(G_2\lor S^7)$ appears in dimension 10, there is a homotopy equivalence $G_2\lor S^7\simeq G_2\times S^7$
% a diagram
% \bse
% \begin{tikzcd}
% G_2\lor S^7 \ar[d,"\p_1'\lor\gamma"]\ar[r,"\alpha\lor\beta"]&[20pt]  \Spin(7)\lor \Spin(7)\ar[r,"\nabla"] \ar[d,"\p_1\lor\p_1"] & \Spin(7)\ar[d]\\
% \Omega^3_0G_2\times\Omega^3S^7\ar[r,"\Omega^3\alpha\times\Omega^3\beta"] &  \Omega^3_0\Spin(7)\times  \Omega^3_0\Spin(7)\ar[r]& \Omega^3_0\Spin(7)
% \end{tikzcd}
% \ese

% \bse
% \begin{tikzcd}
% G_2\lor S^7 \ar[d]\ar[r,"\alpha\times\beta"]&[15pt]  \ar[d,"\simeq"]\Spin(7)\lor \Spin(7)\ar[r,"\nabla"]& \Spin(7)\ar[d]\\
% G_2\times S^7 \ar[d,"\p_1'\times\gamma"]\ar[r,"\alpha\times\beta"]&  \ar[d,"\p_1\times\p_1"]\Spin(7)\times  \Spin(7)\ar[r]& \Spin(7)\ar[d,"\p_1"]\\
% \Omega^3_0G_2\times\Omega^3S^7\ar[r,"\Omega_0^3\alpha\times\Omega^3_0\beta"] &  \Omega^3_0\Spin(7)\times  \Omega^3_0\Spin(7)\ar[r]& \Omega^3_0\Spin(7)
% \end{tikzcd}
% \ese

% Indeed, let $F_k$ be the homotopy fibre of the map $k\gamma\cl S^7\to \Omega S^7$. Then we have a diagram

We therefore have a commutative diagram
\bse
\begin{tikzcd}
G_2\lor S^7\ar[r,"k\p_1'\lor k\gamma"] \ar[d,"\alpha\lor\beta"]&[15pt] \Omega_0^3G_2\times\Omega^3 S^7\ar[d,"\Omega^3\alpha\times\Omega^3\beta","\simeq"']\\
\Spin(7)\ar[r,"k\p_1"]& \Omega^3_0\Spin(7)
\end{tikzcd}
\ese
which induces a map of fibres $\phi\cl M\to\cG_k(\Spin(7))$, where $M$ denotes the homotopy fibre of the map $k\p_1'\lor k\gamma$.

Since the lowest dimensional cell in $G_2\times S^7/(G_2\lor S^7)$ appears in dimension 10, the canonical map $G_2\lor S^7\to G_2\times S^7$ is a homotopy equivalence in dimensions less than 9. It thus follows that $M$ is homotopy equivalent to the homotopy fibre of $k\p_1'\times k\gamma$ in dimensions up to 8. Since the homotopy fibre of $k\p_1'\times k\gamma$ is just the product $\cG_k(G_2)\times F_k$, the composite
% which is commutative upon restriction to the 9-skeleta. Then, in particular, the composite
\bse
C_k\times F_k \longrightarrow \cG_k(G_2)\times F_k \longrightarrow M \xrightarrow{\ \phi\ } \cG_k(\Spin(7))
\ese
induces an isomorphism in mod-2 cohomology in dimensions 1 through 6, and therefore we have
\bse
H^*(\cG_k(\Spin(7)))\cong H^*(C_k)\otimes H^*(F_k), \qquad *\leq 6.
\ese
From the fibration sequence
\bse
%\cdots\longrightarrow
\Omega^4S^7 \longrightarrow F_k \longrightarrow S^7 %\longrightarrow \Omega^3S^7\longrightarrow\cdots
\ese
we see that $H^*(F_k)\cong H^*(\Omega^4S^7)$ in dimensions 1 through 6 for dimensional reasons, and hence we have
\bse
H^*(F_k)\cong \Z/2\Z[y_3,y_6] ,\qquad *\leq 6,
\ese
where $|y_i|=i$, which, in turn, yields
\bse
H^*(\cG_k(\Spin(7))) \cong H^*(C_k)\otimes \Z/2\Z[y_3,y_6], \qquad *\leq 6.
\ese
Since $H^*(F_k)$ does not contribute any generators in degree 4 to $H^*(\cG_k(\Spin(7)))$, the result now follows from Lemma \ref{lem:c_k}. Indeed, the presence of a degree 4 generator allows us to distinguish between the $(4,k)=1$ case and the $2|k$ cases, whereas the vanishing of the Steenrod square $Sq^2$ on the degree 4 generator in $H^*(\cG_k(\Spin(7)))$ coming from  $H^*(C_k)$ can be used to distinguish between the $(4,k)=2$ and $(4,k)=4$ cases.
\eq

\bq[Proof of Theorem \ref{thm:spin7} (b)]
Combine Lemmas \ref{lem:spin7b1} and \ref{lem:spin7b2}.
\eq

\bl
\label{lem:spin8b1}
If $\cG_k(\Spin(8))\simeq\cG_l(\Spin(8))$, then $(4,k)=(4,l)$. 
\el

\bq
As in the proof of Lemma \ref{lem:spin7p=2}, the splitting of $G_2\to\Spin(7)\to S^7$ at the prime 2 implies that there is a 2-local homotopy equivalence
\bse
\mu\circ(\alpha\times\beta)\cl G_2\times S^7 \longrightarrow \Spin(7).
\ese
Since the fibration $\Spin(7)\to \Spin(8)\to S^7$ also splits after localising at any prime, there is a decomposition
\bse
\mu\circ\bigl((\iota\circ\alpha)\times(\iota\circ\beta)\times\gamma\bigr)\cl G_2\times S^7 \times S^7\longrightarrow \Spin(8),
\ese
where $\iota\cl\Spin(7)\to \Spin(8)$ is the inclusion homomorphism and $\gamma$ is a homotopy inverse for the map $\Spin(8)\to S^7$. 

Since the map $\iota\circ \alpha$ is a homomorphism, we have a commutative diagram
\bse
\begin{tikzcd}[column sep=large]
G_2 \ar[d,"\iota\circ\alpha"]\ar[r,"\p_1'"]& \Omega^3_0G_2\ar[d,"\Omega^3(\iota\circ\alpha)"]\\
\Spin(8)\ar[r,"\p_1"] & \Omega^3_0\Spin(8).
\end{tikzcd}
\ese
Furthermore, as $\pi_7(\Omega_0^3G_2)\cong\pi_{10}(G_2)\cong0$, we have
\bse
\pi_7(\Omega^3_0\Spin(8))\cong %\pi_7(\Omega_0^3G_2)\oplus \pi_7(\Omega^3S^7)\cong
\pi_7(\Omega^3S^7)\oplus \pi_7(\Omega^3S^7),
\ese
and thus there are commutative diagrams
\bse
\begin{tikzcd}
S^7 \ar[d,"\iota\circ\beta"]\ar[r,"\delta"]& \Omega^3S^7\times\Omega^3S^7\ar[d,"\Omega^3(\iota\circ\beta)\times\Omega^3\gamma"]\\
\Spin(8)\ar[r,"\p_1"] & \Omega^3_0\Spin(8)
\end{tikzcd}
\quad 
\begin{tikzcd}
S^7 \ar[d,"\gamma"]\ar[r,"\delta'"]& \Omega^3S^7\times\Omega^3S^7\ar[d,"\Omega^3(\iota\circ\beta)\times\Omega^3\gamma"]\\
\Spin(8)\ar[r,"\p_1"] & \Omega^3_0\Spin(8)
\end{tikzcd}
\ese
for some $\delta,\delta'$ representing classes in $\pi_7(\Omega^3S^7\times\Omega^3S^7)\cong(\Z/8\Z)^2$. We therefore have a commutative diagram
\bse
\begin{tikzcd}[column sep=large]
G_2\lor (S^7\lor S^7)\ar[r,"k\p_1'\lor k(\delta \lor \delta')"] \ar[d,"\iota\alpha\lor(\iota\beta\lor\gamma)"]&[17.5pt] \Omega_0^3G_2\times(\Omega^3 S^7\times\Omega^3 S^7)\ar[d,"\Omega^3\iota\alpha\times(\Omega^3\iota\beta\times\Omega^3\gamma)","\simeq"']\\
\Spin(8)\ar[r,"k\p_1"]& \Omega^3_0\Spin(8).
\end{tikzcd}
\ese
Arguing as in the proof of Lemma \ref{lem:spin7b2}, we conclude that 
\bse
H^*(\cG_k(\Spin(7)))\cong H^*(\cG_k(\Spin(8))), \qquad *\leq 6,
\ese
hence the statement follows from Lemma \ref{lem:c_k}.
\eq

\bl
\label{lem:spin8b2}
If $\cG_k(\Spin(8))\simeq\cG_l(\Spin(8))$, then $(7,k)=(7,l)$. 
\el

\bq
%By \cite{theriault-odd}
Localising at $p=7$, we have
\bse
\Spin(8)\simeq \Spin(7)\times S^7\simeq G_2\times S^7\times S^7.
\ese
Applying the functor $\pi_{11}$ and noting that \bse
\pi_{10}(S^7)\cong\pi_{11}(S^7)\cong\pi_{14}(S^7)\cong 0,
\ese
(see, e.g.\ \cite{toda-composition}) we find that the evaluation fibration
\bse
\Spin(8)\xrightarrow{\ \p_k\ } \Omega^3_0\Spin(8)\longrightarrow  \B\cG_k(\Spin(8))\longrightarrow \B \Spin(8)
\ese
reduces to the exact sequence
\bse
\pi_{11}(G_2) \longrightarrow \pi_{11}(\Omega^3_0 G_2)\longrightarrow  \pi_{11}(\B\cG_k(\Spin(8)))\longrightarrow 0.
\ese
Hence the result follows from \cite{g_2}.
\eq

\bl
\label{lem:spin8b3}
If $\cG_k(\Spin(8))\simeq\cG_l(\Spin(8))$ and $k$ and $l$ are multiples of 3, then $(3,k)=(3,l)$.
\el

\bq
By \cite{theriault-odd} or \cite{kishimoto10}, when $k$ is a multiple of 3, there is a 3-local homotopy equivalence
\bse
\cG_k(\Spin(8))\simeq S^7\times\Omega^4S^7\times\cG_k(\Spin(7)).
\ese
Recalling the argument in the proof of Lemma \ref{lem:spin7b1}, we have
\begin{align*}
    \pi_{10}(\cG_k(\Spin(8))) & \cong\pi_{10}(S^7)\oplus\pi_{10}(\Omega^4S^7)\oplus\pi_{10}(\cG_k(\Spin(7)))\\
    & \cong \Z/3\Z\oplus\Z/3\Z\oplus\pi_{10}(\cG_k(\Sp(3)))\\
    & \cong (\Z/3\Z)^2\oplus\Z/(3,k)\Z.
\end{align*}
Hence, if $\cG_k(\Spin(8))\simeq\cG_l(\Spin(8))$ then $\Z/(3,k)\Z\cong \Z/(3,l)\Z$ and thus it must be that $(3,k)=(3,l)$.
\eq

\bq[Proof of Theorem \ref{thm:spin8} (b)]
Combine Lemmas \ref{lem:spin8b1},  \ref{lem:spin8b2} and \ref{lem:spin8b3}.
\eq

\bibliographystyle{plain}
%\addcontentsline{toc}{section}{References}
\bibliography{references}

\begin{bibdiv}
\begin{biblist}

\bib{albanese-spinh}{article}{
      author={Albanese, M.},
      author={Milivojević, A.},
       title={{$\Spin^h$} and further generalisations of spin},
        date={2021},
     journal={J. Geom. Phys.},
      volume={164},
       pages={104174},
}

\bib{atiyah-bott}{article}{
      author={Atiyah, M.~F.},
      author={Bott, R.},
       title={The {Y}ang-{M}ills equations over {R}iemann surfaces},
        date={1983},
     journal={Philos. Trans. Roy. Soc. London Ser. A},
      volume={308},
      number={1505},
       pages={523\ndash 615},
}

\bib{atiyah-bott-shapiro}{article}{
      author={Atiyah, M.~F.},
      author={Bott, R.},
      author={Shapiro, A.},
       title={Clifford modules},
        date={1964},
     journal={Topology},
      volume={3},
      number={1},
       pages={3\ndash 38},
}

\bib{booth-remarks}{article}{
      author={Booth, P.},
      author={Heath, P.},
      author={Morgan, C.},
      author={Piccinini, R.~A.},
       title={Remarks on the homotopy type of groups of gauge transformations},
        date={1981},
     journal={C. R. Math. Acad. Sci. Canada},
      volume={111},
      number={3},
       pages={3\ndash 6},
}

\bib{bryant99}{article}{
      author={Bryant, R.~L.},
      author={Sharpe, E.},
       title={{D}-branes and {$\Spin^c$} structures},
        date={1999},
     journal={Phys. Let. B},
      volume={450},
      number={4},
       pages={353\ndash 357},
}

\bib{claudio-spreafico}{article}{
      author={Claudio, M. H.~A.},
      author={Spreafico, M.},
       title={Homotopy type of gauge groups of quaternionic line bundles over
  spheres},
        date={2009},
     journal={Topol. Its Appl.},
      volume={156},
      number={3},
       pages={643\ndash 651},
}

\bib{cohen-milgram}{inproceedings}{
      author={Cohen, R.~L.},
      author={Milgram, R.~J.},
       title={The homotopy type of gauge theoretic moduli spaces},
        date={1994},
   booktitle={Algebraic topology and its applications},
      editor={Carlsson, Gunnar~E.},
      editor={Cohen, Ralph~L.},
      editor={Hsiang, Wu-Chung},
      editor={Jones, John D.~S.},
      series={Mathematical sciences research institute publications},
   publisher={Springer New York},
       pages={15\ndash 55},
}

\bib{crabb-sutherland}{article}{
      author={Crabb, M.~C.},
      author={Sutherland, W.~A.},
       title={Counting homotopy types of gauge groups},
        date={2000},
     journal={Proc. London Math. Soc.},
      volume={81},
      number={3},
       pages={747\ndash 768},
}

\bib{cutler-sp3}{article}{
      author={Cutler, T.},
       title={The homotopy types of {$\Sp(3)$}-gauge groups},
        date={2018},
     journal={Topol. Its Appl.},
      volume={236},
       pages={44\ndash 58},
}

\bib{cutler-unitary}{article}{
      author={Cutler, T.},
       title={The homotopy types of {$\U(n)$}-gauge groups over {$S^4$} and
  {$\C P^2$}},
        date={2018},
     journal={Homology Homotopy Appl.},
      volume={20},
      number={1},
       pages={5\ndash 36},
}

\bib{cutler-theriault}{article}{
      author={Cutler, T.},
      author={Theriault, S.~D.},
       title={The homotopy types of {$\SU(4)$}-gauge groups},
        date={2019},
     journal={arXiv preprint},
      eprint={1909.04643},
}

\bib{donaldson}{article}{
      author={Donaldson, S.~K.},
       title={An application of gauge theory to four-dimensional topology},
        date={1983},
     journal={J. Differential Geom.},
      volume={18},
      number={2},
       pages={279\ndash 315},
}

\bib{freed-witten}{article}{
      author={Freed, D.~S.},
      author={Witten, E.},
       title={Anomalies in string theory with {D}-branes},
        date={1999},
     journal={Asian J. Math},
      volume={3},
      number={4},
       pages={819\ndash 852},
}

\bib{friedlander}{article}{
      author={Friedlander, E.~M.},
       title={Exceptional isogenies and the classifying spaces of simple {L}ie
  groups},
        date={1975},
     journal={Ann. Math.},
      volume={101},
      number={3},
       pages={510\ndash 520},
}

\bib{gottlieb}{article}{
      author={Gottlieb, D.~H.},
       title={Applications of bundle map theory},
        date={1972},
     journal={Trans. Amer. Math. Soc.},
      volume={171},
       pages={23\ndash 50},
}

\bib{hamanaka-sp2}{article}{
      author={Hamanaka, H.},
      author={Kaji, S.},
      author={Kono, A.},
       title={Samelson products in {$\Sp(2)$}},
        date={2008},
     journal={Topol. Its Appl.},
      volume={155},
      number={11},
       pages={1207\ndash 1212},
}

\bib{hamanaka-kono}{article}{
      author={Hamanaka, H.},
      author={Kono, A.},
       title={Unstable {$K^1$}-group and homotopy type of certain gauge
  groups},
        date={2006},
     journal={Proc. Roy. Soc. Edinburgh Sect. A},
      volume={136},
      number={1},
       pages={149\ndash 155},
}

\bib{hamanaka-s6}{article}{
      author={Hamanaka, H.},
      author={Kono, A.},
       title={Homotopy type of gauge groups of ${\SU}(3)$-bundles over
  ${S}^6$},
        date={2007},
     journal={Topology Appl.},
      volume={154},
      number={7},
       pages={1377\ndash 1380},
}

\bib{konoBSpinc}{article}{
      author={Harada, M.},
      author={Kono, A.},
       title={Cohomology mod 2 of the classifying space of {$\Spin^c(n)$}},
        date={1986},
     journal={Publ. Res. Inst. Math. Sci.},
      volume={22},
      number={3},
       pages={543\ndash 549},
}

\bib{harris61}{article}{
      author={Harris, B.},
       title={On the homotopy groups of the classical groups},
        date={1961},
     journal={Ann. Math},
      volume={74},
      number={2},
       pages={407\ndash 413},
}

\bib{hasui16}{article}{
      author={Hasui, S.},
      author={Kishimoto, D.},
      author={Kono, A.},
      author={Sato, T.},
       title={The homotopy types of {$\PU(3)$}- and {$\PSp(2)$}-gauge groups},
        date={2016},
     journal={Algebr. Geom. Topol.},
      volume={16},
      number={3},
       pages={1813\ndash 1825},
}

\bib{hasui19}{article}{
      author={Hasui, S.},
      author={Kishimoto, D.},
      author={So, T.},
      author={Theriault, S.~D.},
       title={Odd primary homotopy types of the gauge groups of exceptional
  {L}ie groups},
        date={2019},
     journal={Proc. Amer. Math. Soc.},
      volume={147},
      number={4},
       pages={1751\ndash 1762},
}

\bib{huang-highd}{article}{
      author={Huang, R.},
       title={Homotopy of gauge groups over high-dimensional manifolds},
        date={2021},
     journal={Proc. Roy. Soc. Edinburgh Sect. A},
       pages={1\ndash 27},
}

\bib{huang-5d}{article}{
      author={Huang, R.},
       title={Homotopy of gauge groups over non-simply-connected
  five-dimensional manifolds},
        date={2021},
     journal={Sci. China Math.},
      volume={64},
      number={5},
       pages={1061\ndash 1092},
}

\bib{husemoller}{book}{
      author={Husem{\"o}ller, D.},
       title={Fibre bundles},
     edition={3},
      series={Graduate texts in mathematics},
   publisher={Springer-Verlag New York},
        date={1994},
      volume={20},
}

\bib{lang-ehp}{article}{
      author={Jr., G. E.~Lang},
       title={The evaluation map and {EHP} sequences},
        date={1973},
     journal={Pacific J. Math.},
      volume={44},
      number={1},
       pages={201\ndash 210},
}

\bib{jost2008}{book}{
      author={J{\"u}rgen, J.},
       title={Riemannian geometry and geometric analysis},
      series={Universitext},
   publisher={Springer Berlin},
        date={2008},
}

\bib{kishimoto10}{article}{
      author={Kishimoto, D.},
      author={Kono, A.},
       title={Note on mod $p$ decompositions of gauge groups},
        date={2010},
     journal={Proc. Japan Acad. Ser. A},
      volume={86},
      number={1},
       pages={15\ndash 17},
}

\bib{g_2}{article}{
      author={Kishimoto, D.},
      author={Theriault, S.~D.},
      author={Tsutaya, M.},
       title={The homotopy types of {$G_2$}-gauge groups},
        date={2017},
     journal={Topol. Its Appl.},
      volume={228},
       pages={92\ndash 107},
}

\bib{kono91}{article}{
      author={Kono, A.},
       title={A note on the homotopy type of certain gauge groups},
        date={1991},
     journal={Proc. Roy. Soc. Edinburgh Sect. A},
      volume={117},
      number={3-4},
       pages={295\ndash 297},
}

\bib{kono-theriault}{article}{
      author={Kono, A.},
      author={Theriault, S.~D.},
       title={The order of the commutator on {$\SU(3)$} and an application to
  gauge groups},
        date={2013},
     journal={Bull. Belg. Math. Soc. Simon Stevin},
      volume={20},
      number={2},
       pages={359\ndash 370},
}

\bib{kono-tsukuda}{article}{
      author={Kono, A.},
      author={Tsukuda, S.},
       title={A remark on the homotopy type of certain gauge groups},
        date={1996},
     journal={J. Math. Kyoto Univ.},
      volume={36},
      number={1},
       pages={115\ndash 121},
}

\bib{membrillo19}{article}{
      author={Membrillo-Solis, I.~A.},
       title={Homotopy types of gauge groups related to {$S^3$}-bundles over
  {$S^4$}},
        date={2019},
     journal={Topol. Its Appl.},
      volume={255},
       pages={56\ndash 85},
}

\bib{membrillo-theriault}{article}{
      author={Membrillo-Solis, I.~A.},
      author={Theriault, S.~D.},
       title={The homotopy types of {$\U(n)$}-gauge groups over lens spaces},
        date={2021},
     journal={Bol. Soc. Mat. Mex.},
      volume={27},
       pages={40},
}

\bib{mimura67}{article}{
      author={Mimura, M.},
       title={The homotopy groups of {L}ie groups of low rank},
        date={1967},
     journal={J. Math. Kyoto Univ.},
      volume={6},
      number={2},
       pages={131\ndash 176},
}

\bib{mimura95}{incollection}{
      author={Mimura, M.},
       title={Homotopy theory of {L}ie groups},
        date={1995},
   booktitle={Handbook of algebraic topology},
      editor={James, I.~M.},
   publisher={North-Holland},
       pages={951\ndash 991},
}

\bib{mohammadi21}{article}{
      author={Mohammadi, S.},
       title={The homotopy types of {$\PSp(n)$}-gauge groups over {$S^{2m}$}},
        date={2021},
     journal={Topol. its Appl.},
      volume={290},
       pages={107604},
}

\bib{su(4)-s8}{article}{
      author={Mohammadi, S.},
      author={Asadi-Golmankhaneh, M.~A.},
       title={The homotopy types of {$\SU(4)$}-gauge groups over {$S^8$}},
        date={2019},
     journal={Topology Appl.},
      volume={266},
       pages={106845},
}

\bib{su(n)-s6}{article}{
      author={Mohammadi, S.},
      author={Asadi-Golmankhaneh, M.~A.},
       title={The homotopy types of {$\SU(n)$}-gauge groups over {$S^6$}},
        date={2019},
     journal={Topol. Its Appl.},
      volume={270},
       pages={106952},
}

\bib{piccinini1998conjugacy}{book}{
      author={Piccinini, R.~A.},
      author={Spreafico, M.},
       title={Conjugacy classes in gauge groups},
      series={Queen's papers in pure and applied mathematics},
   publisher={Queen's University},
     address={Kingston},
        date={1998},
      volume={111},
}

\bib{rea21}{article}{
      author={Rea, S.},
       title={Homotopy types of gauge groups of {$\PU(p)$}-bundles over
  spheres},
        date={2021},
     journal={J. Homotopy Relat. Struct.},
      volume={16},
       pages={61\ndash 74},
}

\bib{sati12}{article}{
      author={Sati, H.},
       title={Geometry of {$\Spin$} and {$\Spin^c$} structures in the
  {M}-theory partition function},
        date={2012},
     journal={Rev. Math. Phys.},
      volume={24},
      number={3},
       pages={1250005},
}

\bib{so19}{article}{
      author={So, T.},
       title={Homotopy types of gauge groups over non-simply-connected closed
  4-manifolds},
        date={2019},
     journal={Glasgow Math. J.},
      volume={61},
      number={2},
       pages={349\ndash 371},
}

\bib{so19-sun}{article}{
      author={So, T.},
       title={Homotopy types of {$\SU(n)$}-gauge groups over non-spin
  4-manifolds},
        date={2019},
     journal={J. Homotopy Relat. Struct.},
      volume={14},
       pages={787\ndash 811},
}

\bib{sutherland92}{article}{
      author={Sutherland, W.~A.},
       title={Function spaces related to gauge groups},
        date={1992},
     journal={Proc. Roy. Soc. Edinburgh Sect. A},
      volume={121},
      number={1--2},
       pages={185\ndash 190},
}

\bib{theriault-sp2}{article}{
      author={Theriault, S.~D.},
       title={The homotopy types of {$\Sp(2)$}-gauge groups},
        date={2010},
     journal={Kyoto J. Math.},
      volume={50},
      number={3},
       pages={591\ndash 605},
}

\bib{theriault-odd}{article}{
      author={Theriault, S.~D.},
       title={Odd primary homotopy decompositions of gauge groups},
        date={2010},
     journal={Algebr. Geom. Topol.},
      volume={10},
      number={1},
       pages={535\ndash 564},
}

\bib{theriault-moduli1}{article}{
      author={Theriault, S.~D.},
       title={Homotopy decompositions of gauge groups over riemann surfaces and
  applications to moduli spaces},
        date={2011},
     journal={Int. J. Math.},
      volume={22},
      number={12},
       pages={1711\ndash 1719},
}

\bib{theriault-su3}{article}{
      author={Theriault, S.~D.},
       title={The homotopy types of {$\SU(3)$}-gauge groups over simply
  connected 4-manifolds},
        date={2012},
     journal={Publ. Res. Inst. Math. Sci.},
      volume={48},
      number={3},
       pages={543\ndash 563},
}

\bib{theriault-moduli2}{article}{
      author={Theriault, S.~D.},
       title={The homotopy types of gauge groups of nonorientable surfaces and
  applications to moduli spaces},
        date={2013},
     journal={Illinois J. Math.},
      volume={57},
      number={1},
       pages={59\ndash 85},
}

\bib{theriault-su5}{article}{
      author={Theriault, S.~D.},
       title={The homotopy types of {$\SU(5)$}-gauge groups},
        date={2015},
     journal={Osaka J. Math.},
      volume={52},
      number={1},
       pages={15\ndash 31},
}

\bib{theriault-sun}{article}{
      author={Theriault, S.~D.},
       title={Odd primary homotopy types of {$\SU(n)$}-gauge groups},
        date={2017},
     journal={Algebr. Geom. Topol.},
      volume={17},
      number={2},
       pages={1131\ndash 1150},
}

\bib{theriault-so}{article}{
      author={Theriault, S.~D.},
      author={So, T.},
       title={The homotopy types of {$\Sp(2)$}-gauge groups over closed simply
  connected four-manifolds},
        date={2019},
     journal={Proc. Steklov Inst. Math.},
      volume={305},
       pages={287\ndash 304},
}

\bib{toda-composition}{book}{
      author={Toda, H.},
       title={Composition methods in homotopy groups of spheres},
      series={Annals of Mathematics Studies},
   publisher={Princeton University Press},
        date={1962},
      volume={49},
}

\bib{west17}{article}{
      author={West, M.},
       title={Homotopy decompositions of gauge groups over real surfaces},
        date={2017},
     journal={Algebr. Geom. Topol.},
      volume={17},
      number={4},
       pages={2429\ndash 2480},
}

\bib{witten94}{article}{
      author={Witten, E.},
       title={Monopoles and four-manifolds},
        date={1994},
     journal={Math. Res. Lett.},
      volume={1},
      number={6},
       pages={769\ndash 796},
}

\end{biblist}
\end{bibdiv}

\end{document}